\documentclass{amsart}

\usepackage{amssymb}

\newtheorem{thm}{Theorem}%[section]

\newtheorem{prop}[thm]{Proposition}

\theoremstyle{definition}

\newtheorem{say}[thm]{}
\newtheorem{exmp}[thm]{Example}

\newtheorem{note}[thm]{Note}            %\renewcommand{\thenote}{} 
\newtheorem{ack}{Acknowledgments}

\newtheorem{defn-thm}[thm]{Definition--Theorem}  %!!!!!!!!!!!!!!!!!!!!!!!!
\newtheorem{defn-lem}[thm]{Definition--Lemma}  %!!!!!!!!!!!!!!!!!!!!!!!!
\theoremstyle{remark}

\renewcommand{\c}[0]{{\mathbb C}}  

\renewcommand{\o}[0]{{\mathcal O}} 
\newcommand{\p}[0]{{\mathbb P}}

\newcommand{\map}[0]{\dasharrow}
\newcommand{\qtq}[1]{\quad\mbox{#1}\quad}

\begin{document}
\bibliographystyle{amsalpha}

\title[Normal crossing surfaces]{Two examples of surfaces with
 \\ normal crossing singularities}
\author{J\'anos Koll\'ar}

%\today

\maketitle

Let $S$ be a surface over $\c$ with only {\it normal crossing} singularities,
abbreviated as {\it nc}.
That is, each point of $S$ is analytically isomorphic to
one of  3 local models:
smooth point $(x=0)\subset \c^3$, double nc  $(xy=0)\subset \c^3$
or triple nc  $(xyz=0)\subset \c^3$.

The normalization $n:\bar S\to S$ is smooth and the preimage of
the singular locus $D\subset S$ is a nc curve $\bar D\subset \bar S$.
The dualizing sheaf  (or canonical line bundle) $\omega_S$ is locally free
and $n^*\omega_S\cong \omega_{\bar S}(\bar D)$.

The aim of this note is to give examples of nc surfaces
whose canonical line bundle exhibits unexpected behavior.

\begin{prop}\label{prop1}  There is an irreducible, projective,  nc surface
$T_1$ of  general type given in (\ref{exmp1}) whose canonical ring
$$
\sum_{m\geq 0} H^0(T_1, \omega_{T_1}^m)
\qtq{is not finitely generated.}
$$
\end{prop}

\begin{prop}\label{prop2}  There is an irreducible,  nc surface
$T_2$  given in (\ref{exmp2}) such that $ \omega_{T_2}$ is not ample yet
its pull back to the normalization 
$n^* \omega_{T_2}$ is ample.
\end{prop}

The latter answers in the negative
a problem left unresolved in
\cite[III.2.6.2]{ega-11} and posed explicitly  in  \cite[1.12]{vieh}.

I found both of these examples 
while studying the minimal model program for
semi-log-canonical surfaces. The key observation is that
the minimal model program leads to singularities that
satisfy the numerical conditions of log canonicity,
yet no reflexive power of their dualizing sheaf is 
locally free. The pluricanonical forms behave
unexpectedly near such singularities, and this
lies at the heart of both of the examples. 

Semi-log-canonical surfaces naturally appear 
as semi-stable limits of smooth surfaces of general type.
The surface $T_1$ does not arise this way, but, as far as I know, 
there could be examples of nc surfaces which are
smoothable yet whose canonical ring
is not finitely generated.
Indeed, if $g:X\to (c\in C)$ is such a family with nc fiber $X_c$,
then the relative minimal model program produces
$g^m:X^m\to C$ such that $K_{X^m/C}$ is $g^m$-nef,
hence the canonical ring of the central fiber $X^m_c$ 
is  finitely generated. Even if $X\map X^m$ does not contract
any divisor, and hence $X_c$ is birational to $X^m_c$,
the canonical ring of $X_c$ can be different from the 
canonical ring of $X^m_c$. The reason is that flips
in $X\map X^m$ may correspond to blow ups $X_c\leftarrow X^m_c$.
Even for normal log canonical 
surfaces $(S,\Delta)$, the  canonical ring is a birational invariant
only if we declare that all new curves appear in the
boundary $\Delta$  with coefficient 1. As we go from
$X_c$ to $ X^m_c$, the coefficients of the new curves
are dictated by the 3-fold $X$ and are typically less than 1.
Thus all we can assert is that the 
canonical ring of $X^m_c$ is a subring of 
the canonical ring of $X_c$. 

\begin{say}[Gluing along curves]
\label{glue}

Let $S$ be a surface, $C\subset S$ a curve
and $g:C\to C'$ a finite morphism which is
locally analytically a closed immersion.
(Note that this condition holds in the nc case.)

For each $c'\in C'$ glue the local branches of $S$
as dictated by $g:C\to C'$ and finally glue this to
$S\setminus C$.  The resulting surface is denoted by $S/(g)$.
Let $n:S\to S/(g)$ be the corresponding morphism.
If $S$ is normal then $n$ is the normalization of $S/(g)$.

There is no problem doing this as
a complex analytic space, but the (quasi)project\-ivity of  $S/(g)$ can be
quite tricky.  We use the following simple criterion:

{\it Claim.} Assume that $S$ is projective and there is an
ample divisor $H$ such that $H$ intersects $C$  transversally
and $H\cap C=g^{-1}\bigl(g(H\cap C)\bigr)$.
Then $n(H)$ is a Cartier divisor on $S$
and $n(H)$ is ample by \cite[Exrc.III.5.7]{hartsh}.

The existence of $S$ is a very special case of a
 general gluing result \cite{artin}:

Let $X$ be a scheme, $Y\subset X$ a closed subscheme 
and $f:Y\to Y'$ a finite morphism. Then 
there is a unique $F:X\to X'$ such that $F|_Y$ factors through
$f$ and $F$ is maximal with this property.
In general, $X'$ is only an algebraic space.
\end{say}

\begin{say}[Computing sections of $\omega_S^{[m]}$]
\label{sections}

Let $S$ be a reduced surface
and $Z\subset S$ a finite set of points
such that $S\setminus Z$ has only smooth and double nc points. 
As usual, $\omega_S^{[m]}$ denotes the double dual of 
$\omega_S^{\otimes m}$.

Let $n:\bar S\to S$ denote the normalization
and $\bar Z:=n^{-1}(Z)$.
Then $\bar S\setminus \bar Z$ is the normalization
of $S\setminus Z$. Let 
$D\subset S$ be the singular locus and 
$\bar D:=n^{-1}(D)$ its preimage  in $\bar S$.
$\bar D\setminus \bar Z$ is a smooth curve
and 
there is an involution
$\sigma:\bar D\setminus \bar Z\to \bar D\setminus \bar Z$
such that $D\setminus  Z=(\bar D\setminus \bar Z)/\sigma$.

We say that $S$ is obtained from $\bar S$ by the gluing
$\sigma$. Note that $\sigma$ determines
$\bar D\to D$ only on a dense open set.
If we assume in addition that $S$ satisfies Serre's condition $S_2$,
then $\bar S, \bar D$ and $\sigma$ determine $S$ uniquely.

From this description it is easy to compute
the pluricanonical sections:
$$
H^0(S, \omega_S^{[m]})=\bigl\{s\in H^0(\bar S\setminus \bar Z, 
 \omega_{\bar S\setminus \bar Z}^m(m\bar D)) :
\mbox{$s|_{\bar D\setminus \bar Z}$ is 
$\bigl((-1)^m\sigma\bigr)$-invariant}\bigr\}.
\eqno{(\ref{sections}.1)}
$$ 
(See (\ref{loc1}.1) about the sign $(-1)^m$.)
\end{say}

\begin{exmp}\label{exmp1}  Let $A$ be an elliptic curve
with 4 distinct points $p_1, p_2, q_1, q_2\in A$ such that
$p_1+p_2\sim q_1+q_2$.

Let $f:S\to P$ be a genus 2 (irrational ) pencil such that
$\omega_{S}$ is ample  and there are 2 fibers
$F_p\cong A/(p_1\mbox{ identified with } p_2)$ and
 $F_q\cong A/(q_1\mbox{ identified with } q_2)$.

Let $S_1\to S$ be obtained by blowing up
$q_1, q_2\in F_p$ and  $p_1, p_2\in F_q$.
The corresponding exceptional curves are
$E_{q_1}, E_{q_2}, E_{p_1}, E_{p_2}$.

Set  $B:=A/(p_1\mbox{ identified with } p_2, q_1\mbox{ identified with }  q_2)$
and let $C'$ be the union of the 2-nodal curve $B$
 plus 
transversal copies  $\p^1_p, \p^1_q$ through the nodes.

Let $C:=F_p+F_q+E_{p_1}+ E_{p_2}+ E_{q_1}+ E_{q_2}$
and define $g:C\to C'$ to be the identity of $A$ on
$F_p$ and on $F_q$ and isomorphisms
$E_{p_i}\to \p^1_p$ and $E_{q_i}\to \p^1_q$.

Set $T_1:=S_1/(g)$.
$T_1$ has 2 triple points at the 2 nodes of $B$.
\end{exmp}

\begin{exmp}\label{exmp2} Let $C:=(z^2=x^6+2y^6)\subset \p(1,1,3)$
and  $E:=(z^2=xy(x^2+y^2))\subset \p(1,1,2)$
with (hyper)elliptic involutions $\tau_C, \tau_E$.
Let $p\in E$ denote $(0{:}1{:}0)$ and
 $q\in E$ denote $(1{:}0{:}0)$; both are fixed by $\tau_E$.

Set $S:=C\times E/(\tau_C, \tau_E)$. 
Consider the curves
$D_p:=C\times p/(\tau_C, \tau_E)\cong \p^1_{x{:}y}$ and
$D_q:=C\times q/(\tau_C, \tau_E)\cong \p^1_{x{:}y}$.

Let $\sigma:D_p\to D_q$ be the isomorphism which sends
$(x{:}y)$ to $(y{:}x)$
and let $T:=S/(\sigma)$ be the surface obtained by gluing
$D_p$ to $D_q$ using $\sigma$. 
A key property is that $\sigma$ maps nodes to smooth points.

The surface $S$ has 12 ordinary double points, let
$Z\subset T$ be their images. The surface we are looking for is
$T_2:=T\setminus Z$.
\end{exmp}

We start the proofs by the key local computations.
It is then easy to read off the required global properties.

\begin{say}[Local computation 1]\label{loc1}
Let $C_1:=(xy=0)\subset \c^2_{x,y}=:S_1$.
Let $C_{21}:=(u_1=0)\subset \c^2_{u_1,v_1}=:S_{21}$
and $C_{22}:= (v_2=0)\subset  \c^2_{u_2,v_2}=:S_{22}$.
Set $S_2:=S_{21}\coprod S_{22}$ and
$C_2:=C_{21}\coprod C_{22}$

The gluing is defined by $\sigma:C_1\setminus (0,0)\to C_2$
sending $(0,y)\mapsto (0,y)\in C_{21}$ and
$(x, 0)\mapsto  (x,0)\in C_{22}$.

Note that $T:=(S_1\coprod S_2)/\sigma$ is not a nc surface.
Rather, it has a triple point with embedding dimension 4.
A local model is given by 
$$
(t_1=t_2=0)\cup(t_2=t_3=0)\cup(t_3=t_4=0) \subset \c^4.
$$
The isomorphism is given by
$(x,y)\mapsto (0,x,y,0)$,
$(u_1,v_1)\mapsto (v_1,u_1,0,0)$ and
$(u_2,v_2)\mapsto (0,0,v_2,u_2)$.

A local generator of $\omega_{S_{21}}(C_{21})$
is $u_1^{-1}du_1\wedge dv_1$, and the restriction
$\omega_{S_{21}}(C_{21})|_{C_{21}}=\omega_{C_{21}}$
is given by the Poincar\'e residue map
$$
\frac{df}{f}\wedge dg |_{(f=0)} \mapsto dg|_{(f=0)}.
$$
Thus 
$\omega_{S_{21}}^m(mC_{21})|_{C_{21}}=(dv_1)^m
\cdot \o_{C_{21}}$. The situation on $C_{22}$ is similar.

On the other hand, a local generator 
of $\omega_{S_1}(C_1)$ is $(xy)^{-1}dx\wedge dy$.
Its restriction to $C_1$ gives a local generator $\eta$
of $\omega_{C_1}$. 
Note that 
$$
\eta|_{(y=0)}=-\frac{dx}{x}\qtq{and}
\eta|_{(x=0)}=\frac{dy}{y}.\eqno{(\ref{loc1}.1)}
$$
Thus
$$
\omega_{S_1}^m(mC_1)|_{C_1}
=\eta^m\cdot \o_{C_1}.
$$

The interesting feature appears when we compute that
$$
\sigma^* (dv_1)^m =y^m\cdot \eta|_{(x=0)}
\qtq{and}
\sigma^* (du_2)^m =(-x)^m\cdot \eta|_{(y=0)}.
$$
Thus a local section of $\omega_{S_1}^m(mC_1)$
satisfies the gluing condition (\ref{sections}.1) iff
it is contained in
$$
(xy, x^m, y^m)\cdot \bigl(\frac{dx\wedge dy}{xy}\bigr)^m.
$$
Local finite generation fails since the $\o_{S_1}$-algebra
$$
\sum_{m\geq 0}(xy, x^m, y^m)\cdot W^m\subset \c[x,y,W]
\qtq{is not finitely generated,}
$$
where $W$ is a formal variable (or weight) taking care of the grading.
Indeed, for every $m$, the element
$xy\cdot W^m$ needs to be added as a new generator.
\end{say}

\begin{say}[Proof of (\ref{prop1})]
As we noted,
$$
n^*\omega_{T_1}=\omega_{S_1}(F_p+F_q+E_{p_1}+ E_{p_2}+ E_{q_1}+ E_{q_2}),
\qtq{and}
$$
this line bundle has negative degree along the 4 curves
$E_{q_1}, E_{q_2}, E_{p_1}, E_{p_2}$.

That is, the surface
$(S_1, F_p+F_q+E_{p_1}+ E_{p_2}+ E_{q_1}+ E_{q_2})$ is not
log-minimal. Its log-minimal model is $(S, F_p+F_q)$.
Let $T$ be the surface obtained from $S$ by gluing
$F_p$ to $F_q$ by the identity of $A$.
(Note that both $F_p$ and $F_q$ are birational to $A$.)
$T$ is  singular along the 2-nodal curve $B$.
At the nodes $P,Q$  of $B$ we get a singularity as in (\ref{loc1}).

Thus, instead of
thinking of
$\sum_{m\geq 0} H^0(T_1, \omega_{T_1}^m)$ as
a subalgebra of 
$$ 
\sum_{m\geq 0} 
H^0\bigl(S_1, \omega_{S_1}^m(m(F_p+F_q+E_{p_1}+ E_{p_2}+ E_{q_1}+ E_{q_2}
))\bigr),
$$
we  work with
$$
\sum_{m\geq 0} H^0(T, \omega_{T}^{[m]})
=\sum_{m\geq 0} H^0(T_1, \omega_{T_1}^m),
$$
and use the representation
$$
 H^0(T_1, \omega_{T_1}^m)=\bigl\{s\in H^0(S, \omega_{S}^m(m(F_1+F_2)) :
s|_{F_p+F_q}\mbox{ is $\bigl((-1)^m\sigma\bigr)$-invariant}\bigr\}.
$$

Near the two triple points $P,Q\in T$, we are
in the situation described in (\ref{loc1}).
In particular, we know that the $\o_T$-algebra
$\sum_{m\geq 0}  \omega_{T}^{[m]}$ is not
finitely generated, not even locally near $P$ or $Q$.

To go from the local infinite generation
to global infinite generation we 
consider the natural map
$$
\rho:
\sum_{m\geq 0} H^0(T, \omega_{T}^{[m]})\to \sum_{m\geq 0}  \omega_{T}^{[m]}.
$$
Assume that  for all $m\gg 1$ there are
global sections $t_m\in H^0(T, \omega_{T}^{[m]})$
such that
$\rho(t_m)$ is not contained in the
subsheaf of $\omega_{T}^{[m]}$ generated by
the $\omega_{T}^{[i]}$ for $i<m$. Then
$t_m$ is not contained in the
subalgebra generated by
the $H^0(T, \omega_{T}^{[i]})$ for $i<m$, hence
$\sum_{m\geq 0} H^0(T, \omega_{T}^{[m]})$ is not finitely generated.

Since $\omega_S$ is ample and $F_p, F_q$ are nef, we
see that
$\omega_S^m(mF_p+mF_q)(-F_p-F_q)$ is globally generated
for $m\gg 1$. 
Sections of $\omega_S^m(mF_p+mF_q)(-F_p-F_q)$ vanish
along $F_p+F_q$, hence they automatically glue
and descend to sections of $\omega_T^{[m]}$.

Thus if  $s_m\in H^0(S,\omega_S^m(mF_p+mF_q))$ vanishes along
$F_p+F_q$ with multiplicity 1, then 
we obtain a corresponding $t_m\in H^0(T, \omega_{T}^{[m]})$
which, up to a unit, equals $xy\cdot W^m$
in  (\ref{loc1}).
Thus $\sum_{m\geq 0} H^0(T_1, \omega_{T_1}^m)$ is not
finitely generated.

Finally, $T_1$ is projective. To see this note first
that $\omega_{T_1}^{-1}$ is relatively ample
on $T_1\to T$. Thus it is enough to prove that
$T$ is projective.

Note that the pull back of
$\omega_{F_p}$ to $A$ is $\o_A(p_1+p_2)$ and
 the pull back of
$\omega_{F_q}$ to $A$ is $\o_A(q_1+q_2)$.
We assumed that  $p_1+p_2\sim q_1+q_2$, thus there is a divisor
 $H\in |\omega_{S}\otimes f^*(\mbox{very ample})|$
which intersects $F_p$ and $F_q$ transversally in 
points which are interchanged by $\sigma$.  
As noted in (\ref{glue}), $H$ then descends to an ample divisor
on $T$. \qed
\end{say}

\begin{say}[Local computation 2]\label{loc2}
Let $C_1:=(y=0)\subset \c^2_{x,y}=:S_1$.
Let $C_{2}:=(v=w=0)\subset (uv-w^2=0)=:S_{2}\subset \c^2_{u,v,w}$.

The gluing is defined by $\tau:C_2\to C_1$
sending $(u,0,0)\mapsto (u,0)\in C_1$.
As in (\ref{loc1}), we find that
$$
\omega_{S_1}^m(mC_1)|_{C_1}=(dx)^m\cdot \o_{C_1}.
$$
Note that $\omega_{S_2}=u^{-1}du\wedge dw\cdot \o_{S_2}$
and  the equation $(v=0)$ defines $2\cdot C_2$.
Thus
$$
\omega_{S_2}^{2m}(2mC_2)=\bigl(\frac{du\wedge dw}{u}\bigr)^{2m}\cdot 
\frac1{v^m}\cdot \o_{S_2}.
$$
Note that  $w$ vanishes along $C_2$ with multiplicity 1,
hence $w^{-1}du\wedge dw$ restricts to $-du$ on $C_2$.
Thus 
$$
\bigl(\frac{du\wedge dw}{u}\bigr)^{2m}\cdot 
\frac1{v^m}|_{C_2}= 
\bigl(\frac{dw\wedge du}{w}\bigr)^{2m}
\frac{w^{2m}}{u^{2m}v^m}|_{C_2}=
\frac1{u^m}(du)^{2m}.
$$
Hence
$$
\omega_{S_2}^{2m}(2mC_2)|_{C_2}=\frac1{u^m}(du)^{2m}\cdot \o_{C_2}.
$$
Let us now glue $S_1$ to $S_2$ by $\sigma$
to obtain a singular surface $S$ with singular curve $C$.

A section of $\omega_S^{[2m]}$ restricts to a
rational section $h(x)(dx)^{2m}$ of $\omega_C^{2m}$.

Computing on $S_2$ we get that
$h(x)$ seems to be allowed a pole of order $m$
at the origin, but computing on $S_1$ shows that
no pole is allowed. That is
$$
H^0\bigl(S, \omega_S^{[2m]}\bigr)|_C\subset H^0(C, \omega_C^{2m}).
\eqno{(\ref{loc2}.1)}
$$
\end{say}

\begin{say}[Proof of (\ref{prop2})]

Let $n:\bar T_2\to T_2$ denote the normalization.
Then
$$
n^* \omega_{T_2}=\omega_S(D_p+D_q)|_{\bar T_2}.
$$
Note that $\omega_S(D_p+D_q)$ is not Cartier at the
24 nodes, but its reflexive square
$\omega_S^2(2D_p+2D_q)$ is Cartier and ample, thus
$n^* \omega_{T_2}$ is ample.

 $S$ has 6 nodes on $D_p$ and 6 on $D_q$.
Moreover, we chose the isomorphism $\sigma: D_p\cong D_q$
such that every node is matched with a smooth point.
Thus along $D:=n(D_p)=n(D_q)$ we have 12 points whose local models
are as in (\ref{loc2}).

Using  (\ref{loc2}.1) we see that every global section of
$\omega_{T_2}^{2m}$ restricted to $D\cap T_2$ is also 
the restriction of a global section of
$\omega_D^{2m}$ to $D\cap T_2$. Since $D\cong \p^1$,
every global section of
$\omega_{T_2}^m$ vanishes along the double curve $D$.
Thus $\omega_{T_2}$ is not ample. 
In fact,  $\omega_{T_2}$ is not even semi- or weakly-positive in
any sense.

Finally, $T$ is projective since it
has a finite map to $\p^1\times \p^1$ given by
$$
\bigl((x_1{:}y_1{:}z_1),(x_2{:}y_2{:}z_2)\bigr)
\mapsto \bigl((x_1y_1{:}x_1^2+y_1^2), (x_2y_2{:}x_2^2+y_2^2)\bigr).
\qed
$$
\end{say}

\begin{note} The  explicit computation 
in (\ref{loc1}) is a special
case of the following general result:

Let $X$ be a reduced, $S_2$ surface and $F$ a rank 1 sheaf on
$X$. Then the $\o_X$-algebra  $\sum_{m\geq 0} F^{[m]}$ is
finitely generated iff $F^{[m]}$ is locally free for some $m>0$.

It seems that the minimal model of a typical nc surface
has such singularities and its canonical ring is
not finitely generated. 
\end{note}

 \begin{ack} I thank A.\ Corti for a long discussion
concerning an earlier log-surface whose canonical ring
was not finitely generated and B.\ Hassett and B.\ Totaro for correcting 
some errors.
Partial financial support  was provided by  the NSF under grant number 
DMS-0500198. Part of this work was done while the author
visited MSRI and the University of Texas. 
\end{ack}

\bibliography{refs}

\vskip1cm

\noindent Princeton University, Princeton NJ 08544-1000

\begin{verbatim}kollar@math.princeton.edu\end{verbatim}

\end{document}